\documentclass[a4paper, 11pt]{amsart}
\usepackage[utf8]{inputenc}
\usepackage[T1]{fontenc}
\usepackage[english]{babel}
\usepackage[margin=1in]{geometry} 
\usepackage{microtype, xcolor, graphicx}
\usepackage{amsmath, amssymb, amsthm}
\usepackage[linktocpage=true]{hyperref}
\usepackage[capitalise]{cleveref}
\usepackage{enumitem}
\usepackage{tikz-cd} 

\theoremstyle{plain}
\newtheorem{mainthm}{Theorem}

\newtheorem{thm}{Theorem}[section]
\newtheorem{lem}[thm]{Lemma}
\newtheorem{prop}[thm]{Proposition}
\newtheorem{cor}[thm]{Corollary}

\theoremstyle{definition}
\newtheorem{defn}[thm]{Definition}
\newtheorem{exmp}[thm]{Example}

\theoremstyle{remark}
\newtheorem{rem}[thm]{Remark}

\numberwithin{equation}{section}

\DeclareMathOperator{\Fun}{Fun}

\DeclareMathOperator{\Mods}{Mod}
\DeclareMathOperator{\mods}{mod}
\DeclareMathOperator{\Hom}{Hom}
\DeclareMathOperator{\RHom}{\mathbf{R}Hom}
\DeclareMathOperator{\otimesL}{\otimes_R^\mathbf{L}}

\DeclareMathOperator{\Spec}{Spec}

\DeclareMathOperator{\supp}{supp}
\DeclareMathOperator{\id}{id}

\DeclareMathOperator{\hocolim}{\underrightarrow{\mathrm{hocolim}}}

\DeclareMathOperator{\add}{\mathsf{add}}

\DeclareMathOperator{\cosusp}{\mathsf{cosusp}}
\DeclareMathOperator{\aisle}{\mathsf{aisle}}

\DeclareMathOperator{\Coaisle}{Coaisle}
\DeclareMathOperator{\Def}{Def}
\newcommand{\cA}{\mathcal{A}}
\newcommand{\cB}{\mathcal{B}}

\newcommand{\cD}{\mathcal{D}}

\newcommand{\cG}{\mathcal{G}}

\newcommand{\cI}{\mathcal{I}}
\newcommand{\cJ}{\mathcal{J}}
\newcommand{\cK}{\mathcal{K}}

\newcommand{\cM}{\mathcal{M}}

\newcommand{\cS}{\mathcal{S}}
\newcommand{\cT}{\mathcal{T}}
\newcommand{\cU}{\mathcal{U}}
\newcommand{\cV}{\mathcal{V}}
\newcommand{\cW}{\mathcal{W}}
\newcommand{\cX}{\mathcal{X}}

\renewcommand{\a}{\mathfrak{a}} 

\newcommand{\m}{\mathfrak{m}}
\newcommand{\p}{\mathfrak{p}}
\newcommand{\q}{\mathfrak{q}}

\newcommand{\bZ}{\mathbb{Z}}

\newcommand{\bK}{\mathbb{K}} 
\newcommand{\kp}{\kappa(\p)}
\newcommand{\hs}{\mathbf{hs}}

\newcommand{\leftangle}{\left\langle}
\newcommand{\rightangle}{\right\rangle}
\newcommand{\leftcurly}{\left\{}
\newcommand{\rightcurly}{\right\}}

\begin{document}
\title[Telescope conjecture for quiver representations over artinian rings]{Telescope conjecture for quiver representations over artinian rings}

\author{Michal Hrbek}
\address{Michal Hrbek, Institute of Mathematics of the Czech Academy of Sciences, {\v Z}itn{\'a} 25, 115 67 Prague, Czech Republic}
\email{hrbek@math.cas.cz}
\author{Enrico Sabatini}
\address{Enrico Sabatini, Dipartimento di Matematica ``Tullio Levi-Civita'', Universit{\`a} degli Studi di Padova, via Trieste 63, 35121 Padova, Italy}
\email{enrico.sabatini@phd.unipd.it}
\thanks{The first named author was supported by the GA{\v C}R project 23-05148S and the Czech Academy of Sciences (RVO 67985840). The second named author was supported by the Department of Mathematics ``Tullio Levi-Civita'' of the University of Padova.}
\begin{abstract}
    Let $\cD(RC)$ be the derived category of representations of a small category $C$ over a commutative noetherian ring $R$. We study the homotopically smashing t-structures on this category. Specifying our discussion to the stalk categories $\Gamma_\p\cD(RQ)$ for a finite quiver $Q$ and a prime ideal $\p$ of $R$, we prove the telescope conjecture for the derived category of representations of finite quivers over artinian rings. More generally, we prove the same result also outside of the noetherian context, for representations of finite quivers over commutative perfect rings.
\end{abstract}
\maketitle
\setcounter{tocdepth}{1}
\tableofcontents

\section*{Introduction}
The telescope conjecture is a fundamental problem in the study of compactly generated triangulated categories, concerning the relationship between ``large'' and ``small'' objects. Originally formulated in stable homotopy theory by Bousfield and Ravenel \cite{Bou, Rav}, in its classical form it asks whether, for a compactly generated triangulated category, every smashing subcategory is compactly generated. The conjecture has been shown to hold in several algebraic contexts, most notably by Neeman \cite{NeeTC}, who established it for derived categories of commutative noetherian rings. This was later generalized by Antieau and Stevenson \cite{AS} to representations of quivers over commutative noetherian rings, with partial results in the full generality and proving the conjecture in case of Dynkin quivers.

A natural generalization of this problem arises in the setting of t-structures. In analogy with the stable case, one may ask whether every homotopically smashing t-structure is compactly generated. This formulation -- often referred to as the telescope conjecture for t-structures \cite{BHTC, HN} -- recovers the classical conjecture when restricted to stable t-structures. Analogues of the previously mentioned results have been established: Hrbek and Nakamura \cite{HN} proved the conjecture for commutative noetherian rings, and Sabatini \cite{Sab} extended it to representations of Dynkin quivers over such rings. The aim of this paper is to investigate the analogue of the ``partial results'' of \cite{AS} in the general setting of representations of small categories.

Specifically, in \cite{AS} the authors show that there is a bijection between localizing subcategories of $\cD(RC)$ and collections of localizing subcategories of $\cD(\kp C)$ indexed by prime ideals $\p$ in $\Spec(R)$, where $\kp$ is the residue field of $R$ at $\p$ \cite[Corollary 4.3]{AS}. Moreover, for each $\p\in\Spec(R)$ such that the localization $R_\p$ is a regular local ring, this bijection restricts to a bijection between localizing subcategories of the stalk subcategory $\Gamma_\p\cD(RC)$ and localizing subcategories of $\cD(\kp C)$ \cite[Proposition 6.10]{AS}. In this case, for any finite quiver $Q$, the stalk subcategory $\Gamma_\p\cD(RQ)$ satisfies the stable telescope conjecture \cite[Corollary 6.12]{AS}.

In this paper, we establish analogous results for t-structures by showing that there is an injective assignment from homotopically smashing t-structures of $\cD(RC)$ to collections of homotopically smashing t-structures of $\cD(\kp C)$ indexed over $\Spec(R)$.

\begin{mainthm}[\cref{Injection}]\label{ThmA}
    For any commutative noetherian ring $R$ and small category $C$, there is a well-defined injective assignment:
    \begin{equation*}\begin{gathered}
        \Coaisle_{\hs}(\cD(RC))\lhook\joinrel\longrightarrow\prod_{\p\in\Spec(R)}\Coaisle_{\hs}(\cD(\kp C))\\
        \qquad\qquad\cV\longmapsto\left(\add_{\kp C}\leftangle\RHom_R(\kp,\cV)\rightangle\right)_{\p\in\Spec(R)}
    \end{gathered}\end{equation*}
    with its left inverse assigning to a collection of homotopically smashing coaisles $\left(\mathsf{v}(\p)\right)_{\p\in\Spec(R)}$ the smallest cosuspended subcategory of $\cD(RC)$, closed under directed homotopy colimits, that contains $\mathsf{v}(\p)$ for each $\p\in\Spec(R)$.
\end{mainthm}

Moreover, for any finite quiver $Q$ and prime ideal $\p\in\Spec(R)$, such assignment restricts to a bijection between homotopically smashing t-structures of the stalk subcategory $\Gamma_\p\cD(RQ)$ and homotopically smashing t-structures of $\cD(\kp Q)$ (see \cref{Bij}). This proves the telescope conjecture for all stalk subcategories $\Gamma_\p\cD(RQ)$ without assuming regularity of the local ring $R_\p$. Therefore, by trivial gluing techniques, we prove the following:

\begin{mainthm}[\cref{Main}]\label{ThmB}
    For any commutative artinian ring $R$ and finite quiver $Q$, the derived category $\cD(RQ)$ satisfies the telescope conjecture for t-structures. 
\end{mainthm}

Let us stress that, under the regularity assumption of \cite{AS}, the same gluing techniques would only recover the limit case of representations of quivers over fields, for which the result has already been established in \cite{KS} and \cite{AH}. As an application of this result, we exhibit examples of non-hereditary finite dimensional algebras satisfying the telescope conjecture (see \cref{Example}). Moreover, the techniques used to prove \cref{ThmB} can be adapted to establish the telescope conjecture also for representations over certain non-noetherian rings, namely \emph{perfect rings} -- we devote the final section of the paper to this slightly more general result.

\subsection*{Structure of the paper}
The paper is divided into four sections. In \textbf{\cref{Sec1}}, we recall the basics of representations of small categories, algebraic triangulated categories, and t-structures. In \textbf{\cref{Sec2}}, we review the equivalence between homotopically smashing t-structures and those with definable coaisles, and we use the latter notion to prove \cref{ThmA}. In \textbf{\cref{Sec3}}, we specialize to representations of finite quivers and their stalk subcategories, laying the groundwork for our second main result and then we prove \cref{ThmB}. Finally, in \textbf{\cref{Sec4}}, we establish the telescope conjecture for representations of finite quivers over commutative perfect rings.
\subsection*{Acknowledgments}
The authors would like to thank Jorge Vit{\'o}ria for reading a preliminary version of this paper and for providing valuable feedback.

\section{Preliminaries}\label{Sec1}
\subsection{Representations of small categories}
Given a (unitary, associative) ring $R$, let $\Mods(R)$ denote the category of all left $R$-modules. Given a small category $C$, the category $\Mods_R(C)$ of left $C$-modules over $R$ is the category of functors $\Fun(C,\Mods(R))$ or, equivalently, the category of $R$-linear functors $\Fun_R(RC,\Mods(R))$ (see \cite[Lemma 2.7]{AS}), where $RC$ is the $R$-linearization of the small category $C$. Recall from \cite[III 4.2, 4.6, 5.2]{Pop} that $\Mods_R(C)$ is a Grothendieck category with exact products and a generating set of small projective objects. If $\cA$ is an abelian category, let $\cK(\cA)$ denote the homotopy category of cochain complexes over $\cA$ and let $\cD(\cA)$ denote the (unbounded) derived category of $\cA$. We shall use the shorthand notations $\cK(R) = \cK(\Mods(R))$, $\cD(R) = \cD(\Mods(R))$ and $\cK(RC) = \cK(\Mods_R(C))$, $\cD(RC) = \cD(\Mods_R(C))$.

When $R$ is a commutative ring, as it will be from now on, there are two actions of the category $\cD(R)$ on $\cD(RC)$, given by deriving the bifunctors:
\[\_\otimes_R\_:\Mods(R)\times\Mods_R(C)\to\Mods_R(C)\]
\[\Hom_R(\_,\_):\Mods(R)^{op}\times\Mods_R(C)\to\Mods_R(C),\]
where for every $M\in\Mods(R)$, $X\in\Mods_R(C)$ and $c\in C$,
\[(M\otimes_R X)(c)=M\otimes_R X(c)\text{ and }\Hom_R(M,X)(c)=\Hom_R(M,X(c)).\]

\begin{defn}
    For an abelian category $\cA$, a complex $M\in\cK(\cA)$ is called \emph{K-projective} (resp.~\emph{K-injective}) if, for any acyclic object $I\in\cK(\cA)$, it holds that $\Hom_{\cK(\cA)}(M,I)=0$ (resp.~$\Hom_{\cK(\cA)}(I,M)=0$). When $\cA=\Mods(R)$, for a commutative ring $R$, we say that $M\in\cK(R)$ is \emph{K-flat} if the functor $M\otimes_R\_:\cK(R)\to\cK(R)$ preserves acyclicity. Analogously, we say that a complex $Z\in\cK(RC)$ is \emph{$R$-K-flat} (resp.~\emph{$R$-K-injective}) if the functor $\_\otimes_R Z:\cK(R)\to\cK(RC)$ (resp.~$\Hom_R(\_,Z):\cK(R)^\mathrm{op}\to\cK(RC)$) sends acyclic $R$-complexes to acyclic $RC$-complexes.

    Moreover, a K-projective (resp.~K-flat, $R$-K-flat) \emph{resolution} of a given a complex $X\in\cK(RC)$ is a quasi-isomorphism $Z\to X$ in $\cK(RC)$, where $Z$ is K-projective (resp.~K-flat, $R$-K-flat) and we denote any choice of such $Z$ by $\mathbf{p}X$ (resp.~$\mathbf{f}X$, $\mathbf{f_R}X$). Similarly, an $R$-K-injective resolution of $X$ is a quasi-isomorphism $X\to \mathbf{i_R}X$ in $\cK(RC)$, where $\mathbf{i_R}X$ is $R$-K-injective. All of these resolutions always exist for any $X$ by \cite{Spa} and what follows.
\end{defn}

\begin{rem}[{\cite[Lemma 1.7]{Sab}}]\hfill\\
    Define a functor $X\in\Mods_R(C)$ to be \emph{objectwise flat} (resp.~\emph{injective}) if $X(c)$ is a flat (resp.~injective) $R$-module for any $c\in C$. We have that:
    \begin{itemize}[leftmargin=1cm]
        \item[$\boldsymbol{\cdot}$] K-projective complexes of $\cK(RC)$ and bounded above complexes of objectwise flat functors are $R$-K-flat;
        \item[$\boldsymbol{\cdot}$] K-injective complexes of $\cK(RC)$ and bounded below complexes of objectwise injective functors are $R$-K-injective.
    \end{itemize}
\end{rem}

In particular, the following theorem holds:

\begin{thm}[{\cite[Theorem 1.8]{Sab}}]\label{ActionThm}\hfill
\begin{enumerate}
    \item\begin{enumerate}
        \item The bifunctor $\_\otimes_R\_$ admits a left-derived functor $\_\otimesL\_$, which can be computed, independently up to isomorphism, either using K-flat resolutions in the first variable or $R$-K-flat resolutions in the second;
        \item The bifunctor $\Hom_R(\_,\_)$ admits a right-derived functor $\RHom_R(\_,\_)$, which can be computed, independently up to isomorphism, either using K-projective resolutions in the first variable or $R$-K-injective resolutions in the second;
    \end{enumerate}
    \item For any two commutative rings $R$ and $S$ and any $R$-$S$-bicomplex $M$, the functors $M\otimes_S^\mathbf{L}\_$ and $\RHom_R(M,\_)$ form an adjoint pair. In particular, for any $SC$-complex $X$ and $RC$-complex $Y$ there is an isomorphism
    \[\Hom_{\cD(RC)}(M\otimes_S^\mathbf{L}X,Y)\cong\Hom_{\cD(SC)}(X,\RHom_R(M,Y)).\]
\end{enumerate}
\end{thm}

\begin{rem}\label{AssAdj}
    The associativity and adjunction isomorphisms also hold, i.e.~for any $R$-$S$-bicomplex $M$, $R$-complex $L$, $S$-complex $N$, $SC$-complex $X$ and $RC$-complex $Y$ we have the isomorphisms
    \[L\otimesL(M\otimes_S^\mathbf{L}X)\cong(L\otimesL M)\otimes_S^\mathbf{L}X\]
    and
    \[\RHom_R(M\otimes_S^\mathbf{L}N,Y)\cong\RHom_S(N,\RHom_R(M,Y)).\]
    Indeed, this isomorphisms trivially hold for the non-derived functors $\_\otimes\_$ and $\Hom(\_\,\_)$ because they hold in the respective module categories. Then, the isomorphisms above are induced by taking the appropriate resolutions. Explicitly, we have
    \[\mathbf{p}L\otimes_R(M\otimes_S\mathbf{f_S}X)\cong(\mathbf{p}L\otimes_R M)\otimes_S\mathbf{f_S}X\]
    and
    \[\Hom_R(M\otimes_S\mathbf{p}N,\mathbf{i_R}Y)\cong\Hom_S(\mathbf{p}N,\Hom_R(M,\mathbf{i_R}Y)).\]
\end{rem}

For any prime ideal $\p\in\Spec(R)$, we define the adjoint triple $\varphi_\p^\ast\dashv\varphi_\p\dashv\varphi_\p^!$, denoting respectively the base change, the forgetful functor, and the cobase change at $\p$ as follows, via the obvious interpretation of $\kp$ as a bimodule over $R$ and $\kp$:
{\large\[\begin{tikzcd}
\cD(\kp C) \arrow[rrr, "\varphi_\p"] & & & \cD(RC) \arrow[lll, "\varphi_\p^\ast:=\kp \otimesL\_"', bend right] \arrow[lll, "{\varphi_\p^!:=\RHom_R(\kp,\_)}", bend left]
\end{tikzcd}\]}

\noindent By uniqueness of adjoints, the forgetful functor $\varphi_\p$ can be either identified with $\kp\otimes_{\kp}^\mathbf{L}\_$ or $\RHom_{\kp}(\kp,\_)$. In particular, it preserves coproducts and so, by \cite[Theorem 5.1]{Nee}, the functor $\varphi_\p^\ast$ preserves compact objects.

\begin{prop}\label{AdjTrip}
    For any prime ideal $\p$ of $R$, the following hold:
    \begin{enumerate}
        \item $\varphi_\p\circ\varphi_\p^\ast\cong\kp\otimesL\_$ and $\varphi_\p\circ\varphi_\p^!\cong\RHom_R(\kp,\_)$;
        \item If $R$ is noetherian, it holds that
        \[\varphi_\p^\ast\circ\varphi_\p(X)\cong \bigoplus\limits_{i\geq0}X^{n_i}[i] \text{ and } \varphi_\p^!\circ\varphi_\p(X)\cong\prod\limits_{i\geq0}X^{n_i}[-i]\]
        for some $n_i\geq0$ and $n_0=1$;
        \item For any ideal $\q\neq\p$: $\varphi_\q^\ast\circ\varphi_\p=0$ and $\varphi_\q^!\circ\varphi_\p=0$.
    \end{enumerate}
\begin{proof}
    (1) By associativity and adjunction isomorphisms in \cref{AssAdj}, we have that
    \[\varphi_\p\circ\varphi_\p^\ast\cong\left(\kp\otimes_{\kp}^\mathbf{L}\kp\right)\otimesL\_\text{ and }\varphi_\p\circ\varphi_\p^!\cong\RHom_R\left(\kp\otimes_{\kp}^\mathbf{L}\kp,\_\right)\]
    where, trivially, $\kp\otimes_{\kp}^\mathbf{L}{}\kp\cong\kp$.

    \noindent (2) By associativity and adjunction isomorphisms in \cref{AssAdj}, we have that
    \[\varphi_\p^\ast\circ\varphi_\p\cong\left(\kp\otimesL\kp\right)\otimes_{\kp}^\mathbf{L}\_\text{ and }\varphi_\p^!\circ\varphi_\p\cong\RHom_{\kp}\left(\kp\otimesL\kp,\_\right)\]
    where $\kp\otimesL\kp$ is computed by taking a flat resolution $\mathbf{f}\kp$ of $\kp$ over $R$. In fact, we can choose $\mathbf{f}\kp$ to be a resolution of $\kp$ by finite rank free $R_\p$-modules such that $\kp\otimes_R \mathbf{f}\kp$ has zero differentials. Indeed, we can let $\mathbf{f}\kp$ be the minimal projective resolution of the finitely generated $R_\p$-module $\kp$ over the local ring $R_\p$. The rest follows by letting $n_i$ be the dimension of the $(-i)$-th cohomology of $\kp\otimesL\kp$.

    \noindent (3) By \cite[15.1.23]{CFH}, $\kappa(\q)\otimesL\kp=0$.
\end{proof}
\end{prop}

\subsection{Algebraic triangulated categories}
Although, for most of the results in this paper, the generality of a derived category $\cD(\cG)$ of a Grothendieck category $\cG$ -- such as $\cD(RC)$ -- is sufficient to cover our setup and define the \emph{homotopy colimit functor} (see \cite[Example A.2]{HN} for a brief treatment in this context), we will need -- for \cref{MainResSec} -- to define homotopy colimits in compactly generated triangulated categories that are not necessarily derived categories of Grothendieck categories. For this reason, we recall the notion of \emph{algebraic triangulated category}, relying on the equivalence in \cite[Lemma 7.5]{KraAlg}.

\begin{defn}\label{AlgDef}
    A triangulated category $\cT$ is called \emph{algebraic} if there exists a fully faithful exact functor $\cT\longrightarrow\cK(\cA)$ for some additive category $\cA$.
\end{defn}

Note that derived categories of Grothendieck categories are algebraic. Indeed, by the existence of K-injective resolutions, the derived category $\cD(\cG)$ is equivalent to the full triangulated subcategory $\cK_\mathrm{Inj}(\cG)$ of K-injective objects in $\cK(\cG)$.

Moreover, by \cite[Theorem 3.8]{Kel}, any compactly generated algebraic triangulated category is equivalent to the derived category $\cD(\cA)$ of a small dg-category $\cA$, and by \cite[Proposition 1.3.5]{Bec}, $\cD(\cA)$ is equivalent to the homotopy category $\mathrm{Ho}(\mathrm{dgMod}(\cA))$ of the category of dg-modules over $\cA$ endowed with a model structure (we refer to \cite{Kel} for details on dg-categories and to \cite{Hov} for background on model categories). Thus, for any directed small category $I$, we define the direct homotopy colimit functor as the total left derived functor of the direct limit functor:
\[\hocolim_I\:=\mathbf{L}{\underrightarrow{\lim}}_I:\mathrm{Ho}\left(\mathrm{dgMod}(\cA)^I\right)\longrightarrow\mathrm{Ho}(\mathrm{dgMod}(\cA))\cong\cD(\cA)\]

\begin{defn}\label{HocolimDef}
    Given a compactly generated algebraic triangulated category $\cT$, we say that a subcategory $\cX$ is \emph{closed under directed homotopy colimits} if there exists a small dg-category $\cA$ and an equivalence $F:\cT\longrightarrow\cD(\cA)$ such that $F(\cX)$ is closed under homotopy colimits in $\cD(\cA)$.
\end{defn}

It is not a priori clear that this definition is independent of the choice of the dg-category $\cA$, but we will see in \cref{NoAmbRmk} that it is, in the cases relevant to us.

\subsection{t-structures}
Let $\cT$ be a compactly generated algebraic triangulated category. A subcategory of $\cT$ is called \emph{suspended} (resp.~\emph{cosuspended}) if it is closed under direct summands and positive (resp.~negative) shifts and extensions; it is called \emph{cocomplete} (resp.~\emph{complete}) if it is closed under coproducts (resp.~products) and \emph{homotopically smashing} if it is closed under directed homotopy colimits.

\begin{defn}
A \emph{t-structure} in $\cT$ consists of a pair of subcategories $(\cU,\cV)$ satisfying:
\begin{itemize}
    \item[(i)] $\Hom_\cT(\cU,\cV)=0$;
    \item[(ii)] For any $X\in\cT$, there is a distinguished triangle
    \[\begin{tikzcd} U \arrow[r] & X \arrow[r] & V \arrow[r,"+"] & {} \end{tikzcd}\]
    with $U\in\cU$ and $V\in\cV$;
    \item[(iii)] $\cU[1]\subseteq\cU$ (equivalently, $\cV[-1]\subseteq\cV$).
\end{itemize}
\end{defn}

In this case it holds that 
$$\cU={}^\perp\cV := \{X \in \cT \mid \Hom_\cT(X,V) = 0 ~\forall V \in \cV\}$$ 
and 
$$\cV=\cU^\perp := \{X \in \cT \mid \Hom_\cT(U,X) = 0 ~\forall U \in \cU\}.$$ 
We call $\cU$ and $\cV$ the \emph{aisle} and the \emph{coaisle} of the t-structure, respectively. In particular, any aisle is a cocomplete suspended subcategory closed under homotopy colimits (see \cite[Proposition 5.2]{SSV}) and, by \cite{KV}, a suspended subcategory of $\cT$ is an aisle if and only if the inclusion functor has a right adjoint, which we call the (\emph{left}) \emph{truncation functor} of the t-structure; the dual holds for coaisles and cosuspended subcategories.

\begin{defn}
Given a set of objects $\cS\subseteq\cT$ we will denote by $\aisle_\cT\leftangle\cS\rightangle$ the smallest aisle of $\cT$ containing the object in $\cS$, which exists by \cite[Theorem 2.3]{NeeSet}. Then, a t-structure $(\cU,\cV)$ in $\cT$ is called:
\begin{itemize}[label=$\boldsymbol{\cdot}$]
    \item \emph{Compactly generated} if $\cU=\aisle_\cT\leftangle\cS\rightangle$, for some $\cS\subseteq\cT^c$;
    \item \emph{Homotopically smashing} if $\cV$ is closed under directed homotopy colimits;
    \item \emph{Stable} if $\cU$ and $\cV$ are triangulated subcategories of $\cT$, equivalently, if $\cU$ is closed under $[-1]$ or $\cV$ is closed under $[1]$.
\end{itemize}
By \cite[Proposition 7.2]{SSV}, any compactly generated t-structure is homotopically smashing, and we say that the (\emph{stable}) \emph{telescope conjecture} holds in $\cT$ if any (stable) homotopically smashing t-structure is compactly generated.
\end{defn}

Note that strict localizing subcategories, in the sense of Bousfield and Krause, are precisely the aisles of stable t-structures and, by \cite[Theorem A]{KraSma}, a stable t-structure $(\cU,\cV)$ is homotopically smashing if and only if $\cU$ is a smashing subcategory. Thus, our definition of the stable telescope conjecture recovers the classical one. 

\section{Homotopically smashing t-structures}\label{Sec2}
\subsection{Definable coaisles}
Under mild assumptions, the coaisles of homotopically smashing t-structures turn out to be definable in the following sense.

\begin{defn}
    A subcategory of a compactly generated triangulated category $\cT$ is called \emph{definable} if it is of the form
    \[\cI^\perp:=\{X\in\cT\mid\Hom_\cT(f,X)=0\text{ for any }f\in\cI\}\]
    for a set of maps between compact objects $\cI\subseteq\mathrm{Mor} (\cT^c)$.
\end{defn}

Note that definable subcategories are not assumed a priori to satisfy any closure conditions such as closure under shifts or extensions. However, these can be imposed by requiring corresponding properties on the set $\cI$. In particular, by \cite[Proposition 8.17]{SS}, definable subcategories are in bijection with the so-called \emph{saturated ideals} $\cI$ of $\cT^c$. Moreover, by the arguments in \cite[Corollary 8.19]{SS}, \emph{idempotent} saturated ideals give rise to extension-closed definable subcategories (although it is not yet known whether this correspondence is bijective). Finally, closure under shifts can also be recovered: indeed, by \cite[Theorem 8.16]{SS} (resp.~\cite[Corollary 12.5]{Kra}), cosuspended (resp.~triangulated) definable subcategories are in bijection with \emph{suspended} (resp.~\emph{exact}) ideals. In the latter case, the definable subcategories fit in the right-hand sides of (stable) t-structures.

The following shows that in our setting these t-structures are precisely the homotopically smashing ones. Because of this property, it may happen that we will refer to homotopically smashing t-structures also as t-structures with \emph{definable coaisles}.

\begin{prop}[{\cite{SS, LV}}]\label{EqProp}
Let $\cT$ be a compactly generated algebraic triangulated category, then a subcategory $\cV\subseteq\cT$ is the coaisle of a homotopically smashing t-structure if and only if it is a cosuspended definable subcategory.
\begin{proof}
    By \cite[Remark 8.9, (3)$\Leftrightarrow$(2')]{SS}, we have that any homotopically smashing t-structure has a definable coaisle, which by hypothesis is also cosuspended. Moreover, by \cite[Theorem 4.7, Proposition 5.10]{LV}, any cosuspended definable subcategory is the coaisle of a t-structure and it is closed under directed homotopy colimits.
\end{proof}
\end{prop}

Next remark clarifies why we should not care about the ambiguity of \cref{HocolimDef} in the following.

\begin{rem}\label{NoAmbRmk}
    By \cref{EqProp}, a subcategory $\cV\subseteq\cT$ of a compactly generated algebraic triangulated category is a homotopically smashing coaisle if and only if it is the right orthogonal to a suspended ideal $\cI\subseteq\cT^c$. In this case, the property of being closed under directed homotopy colimits does not depend on the choice of a dg-enhancement of $\cT$, as it only depends on the full subcategory of compact objects $\cT^c$. For the same reason, when studying cosuspended subcategories closed under directed homotopy colimits, the choice of dg-enhancement is irrelevant, provided the subcategory is also definable.

    This issue becomes even more moot in the case of a compactly generated derived category of a Grothendieck category -- such as $\cD(RC)$ -- since in this case, by \cite{CS}, the choice of the dg-enhancement is actually unique.
\end{rem}

In view of \cref{EqProp}, the (stable) telescope conjecture can be reformulated in terms of suspended (resp.~exact) ideals (we refer to \cite[Theorem 13.4]{Kra} for the stable case).

\begin{prop}\label{TCideal}
    Let $\cT$ be a compactly generated algebraic triangulated category. Then the following are equivalent:
    \begin{enumerate}
        \item Every homotopically smashing t-structure of $\cT$ is compactly generated;
        \item Every suspended ideal in $\cT^c$ is generated by identity maps.
    \end{enumerate} 
\begin{proof}
    By \cite[Theorem 8.16]{SS}, any definable coaisle is right orthogonal to a set of compact objects if and only if for any suspended ideal $\cI\subseteq\cT^c$ there exists a set of compact objects $\cS\subseteq\cT^c$ such that $\cI^\perp=\cS^\perp$. In this case, any map in $\cI$ factors through an object of $\aisle_\cT\langle\cS\rangle$ and thus, by \cite[Proposition 8.27]{SS}, through an object of $\aisle_\cT\langle\cS\rangle\cap\cT^c$. Therefore, the identities on these objects generate $\cI$ . On the other hand, if $\cI$ is generated by identity maps we can take $\cS$ to be the set of compacts such that $\{\id_S\}_{S\in\cS}$ generates $\cI$.
\end{proof}
\end{prop}

\subsection{Stalk subcategories}
From now on we assume $R$ to be a commutative noetherian ring, unless stated otherwise. We recall that the suitable versions of the \emph{local-to-global principle} and the \emph{minimality of stalk subcategories} hold in $\cD(RC)$.

\begin{defn}
    For any prime ideal $\p\in\Spec(R)$, we let $\Gamma_\p R$ be the tensor-idempotent residue object at $\p$ in $\cD(R)$ (we refer to \cite[Section 3.4]{Stour} for its definition). Explicitly, $\Gamma_\p R$ is the local cohomology complex $\mathbf{R}\Gamma_{V(\p)}R_\p$ of the local ring $R_\p$. For any small category $C$, we call the \emph{stalk subcategory} of $\cD(RC)$ at $\p$ the essential image $\Gamma_\p\cD(RC):=\Gamma_\p R\otimes_R^\mathbf{L} \cD(RC)$ of the endofunctor $\Gamma_\p R\otimesL\_$ on $\cD(RC)$.
\end{defn}

\begin{rem}\label{SuppProp}
    By \cite[Remark 1.12, Proposition 1.13]{Sab}, defining the support of an object $X\in\cD(RC)$ as $\supp_{RC}(X):=\{\p\in\Spec(R)\mid\kp\otimesL X\neq0\}$, we have that for any $\p\in\Spec(R)$
    \[X\in\Gamma_\p\cD(RC) \text{ if and only if } \supp_{RC}(X)\subseteq\{\p\}.\]
    Namely, the stalk subcategory of $\cD(RC)$ at $\p$ is the full subcategory of $\cD(RC)$ containing the objects supported on $\{\p\}$.
\end{rem}

Given a collection of objects $\cS\subseteq\cT$ we will denote by $\cosusp_\cT\leftangle\cS\rightangle$ the smallest cosuspended subcategory of $\cT$ containing the object in $\cS$ and we will use the decorations $\Pi,\hs$ in the superscript to mean that the cosuspended subcategory is the smallest among the complete and/or the homotopically smashing, respectively.

The ``coaisle versions'' of the \emph{local-to-global principle} and the \emph{minimality of stalk subcategories} in $\cD(RC)$ are stated as follows.

\begin{prop}[{\cite[Theorem 2.2, Remark 2.4]{Sab}}]\label{LTGM}\hfill
\begin{enumerate}
    \item For any complex $Y\in\cD(RC)$:
    \[\cosusp_{RC}^\hs\leftangle Y\rightangle=\cosusp_{RC}^\hs\leftangle\Gamma_\p Y\mid\p\in\Spec(R)\rightangle\]
    \item For any complex $Y\in\Gamma_\p\cD(RC)$:
    \[\cosusp_{RC}^\hs\leftangle Y\rightangle\subseteq\cosusp_{RC}^\hs\leftangle\RHom_R(\kp, Y)\rightangle\]
\end{enumerate}
\end{prop}

\subsection{}
Let us now set the background for the study of homotopically smashing t-structures over $\cD(RC)$ by considering their ``local pieces'' over $\cD(\kp C)$, for any prime ideal $\p\in\Spec(R)$. In particular, we will show that any definable coaisle of $\cD(RC)$ can be decomposed into a collection of definable coaisles of $\cD(\kp C)$ indexed over $\Spec(R)$ and that this assignment is injective.

Denote the disjoint union of the collections of definable coaisles of $\cD(\kp C)$ by
\[\mathfrak{C}_{\Def}=\bigsqcup_{\p\in\Spec(R)}\Coaisle_{\Def}(\cD(\kp C)).\]
By a map $\mathsf{f}:\Spec(R)\longrightarrow\mathfrak{C}_{\Def}$ we will always implicitly mean that $\mathsf{f}(\p)\in\Coaisle_{\Def}(\cD(\kp C))$, i.e.~we will only consider sections of the map $\pi:\mathfrak{C}_{\Def}\to\Spec(R)$ assigning to a definable coaisle of $\cD(\kp C)$ the prime $\p$. Note that, by elementary set theory, considering these maps is equivalent to considering collections of definable coaisles of $\cD(\kp C)$ indexed over $\Spec(R)$ (as in \cref{ThmA}); however, we adopt the former setting to remain consistent with \cite{AS, Sab}.

\begin{thm}\label{Injection}\hfill
    \begin{enumerate}
        \item For any definable coaisle $\cV=\cI^\perp$ of $\cD(RC)$ and any prime ideal $\p\in\Spec(R)$, it holds that $\add_{\kp C}\leftangle\varphi_\p^!\cV\rightangle=\left(\varphi_\p^\ast\cI\right)^\perp$. Furthermore, setting $F(\cV)(\p) := \add_{\kp C}\leftangle\varphi_\p^!\cV\rightangle$ gives rise to a well-defined assignment
        \[F:\Coaisle_{\Def}(\cD(RC))\longrightarrow\{\Spec(R)\to\mathfrak{C}_{\Def}\}.\]
        \item For a map $\mathsf{v}:\Spec(R)\longrightarrow\mathfrak{C}_{\Def}$, let $G(\mathsf{v}):=\cosusp_{RC}^\hs\langle\varphi_\p\mathsf{v}(\p)\mid\p\in\Spec(R)\rangle$. Then, for any definable coaisle $\cV\subseteq\cD(RC)$ we have $G\circ F(\cV)=\cV$. In particular, $F$ is injective.
    \end{enumerate}
\begin{proof}
    (1) Let us first prove the equality. For any morphism in $\cI$, say $f:A\to B$, and $Y\in\cV$, by adjunction, we have the following commutative square:
    \[\begin{tikzcd}
        {\Hom_{\cD(\kp C)}\left(\varphi_\p^\ast B,\varphi_\p^! Y\right)} \arrow[dd, "\|\wr", phantom] \arrow[rr, "{\Hom(\varphi_\p^* f, \varphi_\p^! Y)}"] & {} & {\Hom_{\cD(\kp C)}\left(\varphi_\p^\ast A,\varphi_\p^! Y\right)} \arrow[dd, "\|\wr", phantom] \\
         & & \\
        {\Hom_{\cD(RC)}\left(B,\varphi_\p \varphi_\p^! Y\right)} \arrow[rr, "{\Hom(f, \varphi_\p \varphi_\p^! Y)}"] & {} & {\Hom_{\cD(RC)}\left(A,\varphi_\p \varphi_\p^! Y\right)}
    \end{tikzcd}\]
    In particular, $\Hom_{\cD(\kp C)}\left(\varphi_\p^\ast\cI,\varphi_\p^!\cV\right)=0$ if and only if $\Hom_{\cD(RC)}\left(\cI,\varphi_\p \varphi_\p^!\cV\right)=0$. Then, since by \cref{AdjTrip} (1) we have $\varphi_\p\varphi_\p^!\cV\cong\RHom_R(\kp,\cV)$ and, by \cite[Proposition 1.11 (2)]{Sab} we have that $\RHom_R(\kp,\cV)$ is contained in $\cV$, it follows that $\Hom_{\cD(\kp C)}\left(\varphi_\p^\ast\cI,\varphi_\p^!\cV\right)=0$. Thus, since $\left(\varphi_\p^\ast\cI\right)^\perp$ is closed under summands, we get that $\add_{\kp C}\leftangle\varphi_\p^!\cV\rightangle\subseteq\left(\varphi_\p^\ast\cI\right)^\perp$. For the other containment, let $Y\in\left(\varphi_\p^\ast\cI\right)^\perp$. By adjunction, we have that $\varphi_\p Y\in\cI^\perp=\cV$, thus $\varphi_\p^!\varphi_\p Y$ lies in $\varphi_\p^!\cV$ and so, by \cref{AdjTrip} (2), $Y$ belongs to $\add_{\kp C}\leftangle\varphi_\p^!\cV\rightangle$. Now let us prove that the assignment is well-defined, i.e.~that $F(\cV)(\p)$ is a definable coaisle of $\cD(\kp C)$. Since it is evidently definable and closed under summands and negative shifts, it is sufficient to check the closure under extensions. Without loss of generality we consider a triangle
    \[\begin{tikzcd} \varphi_\p^! X \arrow[r] & Y \arrow[r] & \varphi_\p^! Z \arrow[r, "+"] & {} \end{tikzcd}\]
    with $X,Z\in\cV$ and, applying the forgetful functor $\varphi_\p$, as above, we get that $\varphi_\p Y\in\cV$ and so $Y$ is in $\add_{\kp C}\leftangle\varphi_\p^!\cV\rightangle$.
    
    \noindent (2) Notice that, for any definable coaisle $\cV$ of $\cD(RQ)$ we have
    \[G\circ F(\cV)=\cosusp_{RQ}^\hs\leftangle\varphi_\p\varphi_\p^!\cV\mid\p\in\Spec(R)\rightangle\]
    Thus, it is contained in $\cV$. On the other hand, for any complex $Y$ in $\cV$ and $\p\in\Spec(R)$, by \cite[Proposition 1.11 (3)]{Sab}, the complex $\Gamma_\p Y$ is also in $\cV$, thus $\RHom_R(\kp,\Gamma_\p Y)=\varphi_\p\varphi_\p^!\Gamma_\p Y$ is contained in $G\circ F(\cV)$ and, by minimality and local-to-global principles (\cref{LTGM}), $Y\in G\circ F(\cV)$.
\end{proof}
\end{thm}

\begin{rem}
    Actually, one can prove that the injection $F$ is just a ``shadow'' of two different, more general, Galois connections. Indeed, defining
    \[\mathfrak{C}=\bigsqcup_{\p\in\Spec(R)}\mathrm{Cosusp}(\cD(\kp C)) \text{ and } \mathfrak{D}=\bigsqcup_{\p\in\Spec(R)}\Def(\cD(\kp C)),\]
    where $\Def(\cT)$ denotes the set of all definable subcategories of $\cT$, the following assignments form two Galois connections
    \[\begin{tabular}{c c c c}
        $\mathrm{Cosusp}(\cD(RC))\longleftrightarrow\{\Spec(R)\longrightarrow\mathfrak{C}\}$ & & & $\Def(\cD(RC))\longleftrightarrow\{\Spec(R)\longrightarrow\mathfrak{D}\}$ \\
        \begin{tikzcd} \cB \arrow[r, "F^\mathfrak{C}", maps to] & \left(\p\mapsto\cosusp_{\kp C}\leftangle\varphi_\p^!\cB\rightangle\right) \end{tikzcd} & & & \begin{tikzcd} \leftcurly f\in\cD^c(RC)\mid\varphi_\p^\ast f\in\mathsf{i}(\p)\rightcurly^\perp & \mathsf{i}^\perp \arrow[l, "G^\mathfrak{D}"', maps to] \end{tikzcd} \\
        \begin{tikzcd} \leftcurly X\in\cD(RC)\mid\varphi_\p^! X\in\mathsf{b}(\p)\rightcurly & \mathsf{b} \arrow[l, "G^\mathfrak{C}"', maps to] \end{tikzcd} & & & \begin{tikzcd} \cI^\perp \arrow[r, "F^\mathfrak{D}", maps to] & \left(\p\mapsto\left(\varphi_\p^\ast\cI\right)^\perp\right) \end{tikzcd}
    \end{tabular}\]
    In particular:
    \begin{itemize}
        \item[$\boldsymbol{\cdot}$] For any cosuspended subcategory $\cB\subseteq\cD(RC)$ and any map $\mathsf{b}:\Spec(R)\to\mathfrak{C}$, it holds that $F^\mathfrak{C}(\cB)\subseteq \mathsf{b}$ if and only if $\cB\subseteq G^\mathfrak{C}(\mathsf{b})$;
        \item[$\boldsymbol{\cdot}$] For any definable subcategory $\cI^\perp\subseteq\cD^c(RC)$ and any map $\mathsf{i}^\perp:\Spec(R)\to\mathfrak{D}$, it holds that $G^\mathfrak{D}(\mathsf{i}^\perp)\subseteq\cI^\perp$ if and only if $\mathsf{i}^\perp\subseteq F^\mathfrak{D}(\cI^\perp)$.
    \end{itemize}
    What \cref{Injection} shows is that the a priori distinct assignments $F^\mathfrak{C}$ and $F^\mathfrak{D}$ coincide on subcategories that are both cosuspended and definable, and in this case, they share the same left inverse $G$.
\end{rem}

\section{Representation of quivers}\label{Sec3}
When the small category $C$ is the free category on a finite quiver $Q$, the category $\Mods_R(C)$ is equivalent to the category of modules $\Mods(RQ)$, where $RQ$ is the free $R$-algebra on the set of paths of $Q$ with the composition of paths as internal product. In this case, the functors $X:C\to\Mods(R)$ can be seen as representations of the quiver $Q=(Q_0,Q_1)$ over the ring $R$, i.e.~$X=(X_i,X_\alpha)_{i\in Q_0,\alpha\in Q_1}$ where $X_i$ is an $R$-module for any vertex $i\in Q_0$ and $X_\alpha:X_i\to X_j$ is an $R$-linear map for any arrow $\alpha:i\to j\in Q_1$. Recall that, by \cite[Theorem 4.1.4]{Ben}, under the assumption that $Q$ is finite, for any field $\bK$ the path algebra $\bK Q$ is hereditary.

The aim of this section is to prove that, in this setting, for any prime ideal $\p\in\Spec(R)$, the assignment of \cref{Injection} restricts to a bijection between the class of definable coaisles of the stalk subcategory $\Gamma_\p\cD(RQ)$ and the class of definable coaisles of $\cD(\kp Q)$. Moreover, in virtue of these bijections, the former category satisfies the telescope conjecture.

\begin{rem}\label{Gammap}
    By \cref{AlgDef}, it is clear that full subcategories of algebraic triangulated categories are still algebraic. Thus, by \cite[Lemma 6.9]{AS}, for any prime ideal $\p\in\Spec(R)$, the stalk subcategory $\Gamma_\p\cD(RQ)$ is a compactly generated algebraic triangulated categories such that
    \[(\Gamma_\p\cD(RQ))^c=\Gamma_\p\cD(RQ)\cap\cD^c(R_\p Q).\]
\end{rem}

\subsection{Lift of compacts}
First, we want to show a construction for ``lifting'' a compact object $S\in\cD(\kp Q)$ to a compact object $\widetilde{S}\in\Gamma_\p\cD(RQ)$, in such a way that its ``restriction'' $\kp\otimes_R\widetilde{S}$ is an object closely related to the original object $S$.

Before we start, we need to recall the definition of a \emph{Koszul complex} and some of its properties.

\begin{defn}
    Let $\a$ be an ideal of $R$ and $\{p_1,\ldots,p_n\}$ a set of generators for $\a$. Then, the \emph{Koszul complex} (on $R$) at $\a$ is defined as
    \[K(\a)=\bigotimes_{i=1}^n R\xrightarrow{p_i\cdot}R\]
    where the complex $R\xrightarrow{p_i\cdot}R$ is concentrated in degrees $-1$ and $0$ and the map is the multiplication by $p_i$. 
\end{defn}

\begin{rem}\label{Koszul}\hfill
\begin{itemize}
    \item[$\boldsymbol{\cdot}$] Although, in some situations, the definition may be independent of the choice of generators (see \cite[14.4.29]{CFH}), in general this is not the case, so for each ideal $\p$ we choose a set of generators of minimal cardinality;
    \item[$\boldsymbol{\cdot}$] For any prime ideal $\q$, it holds that 
    \[\kappa(\q)\otimesL K(\p)=\begin{cases} \bigoplus\limits_{i=0}^n\kappa(\q)^{\binom{n}{i}}[i] & \text{if } \p\subseteq\q \\
    \quad\quad 0 & \text{otherwise} \end{cases}\]
    Indeed, by \cite[15.1.10]{CFH} the complex is non-zero if and only if $\p\subseteq\q$ and in this case, tensoring by $\kappa(\q)$ annihilates multiplications by elements of $\p$, so it annihilates all the differentials of $K(\p)$.
\end{itemize}
\end{rem}

We are now ready to provide our lifting construction.

\begin{prop}\label{Lift}
    Let $\p\in\Spec(R)$. Then, for any compact object $S\in\cD(\kp Q)$ there exists:
    \begin{enumerate}
        \item $\widehat{S}\in\cD^c(R_\p Q)$ such that $\kp\otimes_R\widehat{S}\cong S$;
        \item $\widetilde{S}\in(\Gamma_\p\cD(RQ))^c$ such that $\kp\otimes_R\widetilde{S}\in\add_{\kp Q}\langle S[\geq0]\rangle$.
    \end{enumerate}
\begin{proof}
    Since the path algebra $\kp Q$ is hereditary, any complex $X\in\cD(\kp Q)$ is isomorphic to $\bigoplus_{i\in\bZ}H^i(X)[-i]$, see \cite[\S 1.6]{KraAlg}. So instead of starting with a complex $S\in\cD^c(\kp Q)$, without loss of generality, we can assume $S\in\mods(\kp Q)$ to be a finitely presented module and then extend the construction to complexes by direct sums and shifts.

    \noindent(1) For such $S\in\mods(\kp Q)$, we have an exact sequence
    \[\begin{tikzcd} 0 \arrow[r] & P \arrow[r, "A"] & {\kp Q}^n \arrow[r] & S \arrow[r] & 0 \end{tikzcd}\]
    where $P=(\kp^{m_i},P_\alpha)_{i\in Q_0,\,\alpha\in Q_1}$ is a finitely generated projective module. In particular, this exact sequence induces a distinguished triangle
    \[\begin{tikzcd} P \arrow[r, "A"] & {\kp Q}^n \arrow[r] & S \arrow[r, "+"] & {} \end{tikzcd}\]
    By the Krull-Schmidt property, $P$ is isomorphic to a direct sum of indecomposable projective modules which, by \cite{Rin}, admit a presentation by $0,1$-matrices on the arrows. In particular, the entries of $P_\alpha$ involve just zero and identity maps for any $\alpha\in Q_1$, thus $P$ can be lifted to a finitely generated projective $R_\p Q$-module $\widehat{P}=(R_p^{m_i},P_\alpha)_{i\in Q_0,\,\alpha\in Q_1}$. By \cite[Theorem 4.1]{CB}, the morphism $A$, given by a collection of matrices $(A_i)_{i\in Q_0}$ with entries in $\kp$, lifts to a morphism $\widehat{A}:\widehat{P}\to{R_\p Q}^n$, given by a collection of matrices $(\widehat{A}_i)_{i\in Q_0}$ with entries in $R_\p$ such that $\kp\otimes_R\widehat{A}=A$. We can complete $\widehat{A}$ to a triangle in $\cD^c(R_\p Q)$
    \begin{equation*}\label{star}\tag{$\star$}
        \begin{tikzcd} \widehat{P} \arrow[r, "\widehat{A}"] & {R_\p Q}^n \arrow[r] & \widehat{S} \arrow[r, "+"] & {} \end{tikzcd}
    \end{equation*}
    and, by construction, we get that $\kp\otimes_R\widehat{S}$ is quasi-isomorphic to $S$.

    \noindent(2) Let $K(\p)$ be the Koszul complex at $\p$. Since it is compact in $\cD(R)$, the complexes $K(\p)\otimes_R\widehat{P}$ and $K(\p)\otimes_R{R_\p Q}^n$ are compact in $\cD(R_\p Q)$ and, by \cref{Koszul}, are supported on $\{\p\}$, thus they belong to $\Gamma_\p\cD(R_\p Q)$ (see \cite[Proposition 1.13 (3)]{Sab}). By \cref{Gammap}, it follows that both the objects are compact in $\Gamma_\p\cD(R_\p Q)$. Define $\widetilde{S}:=K(\p)\otimes_R\widehat{S}$ and note that, by tensoring the triangle (\ref{star}) by $K(\p)$, we get that $\widetilde{S}$ lies in $(\Gamma_\p\cD(R_\p Q))^c$. Moreover, by \cref{Koszul}, we have that
    \[\kp\otimes_R\widetilde{S}=\bigoplus\limits_{i=0}^n\kp^{\binom{n}{i}}[i]\otimes_R\widehat{S}=\bigoplus\limits_{i=0}^n S^{\binom{n}{i}}[i].\]
    In particular, $\kp\otimes_R\widetilde{S}\in\add_{\kp Q}\langle S[\geq0]\rangle$.
\end{proof}
\end{prop}

\subsection{Main results}\label{MainResSec}
The adjoint triple $\varphi_\p^\ast\dashv\varphi_\p\dashv\varphi_\p^!$ introduced in \cref{Sec1} restricts to an adjunction
\[\begin{tikzcd}
\cD(\kp Q) \arrow[rr, "\overline{\varphi}_\p"] & & \Gamma_\p\cD(RQ) \arrow[ll, "\varphi_\p^\ast\circ\,\iota_\p"', bend right] \arrow[ll, "\varphi_\p^!\circ\,\iota_\p", bend left]
\end{tikzcd}\]
where $\iota_\p:\Gamma_\p\cD(RQ)\hookrightarrow\cD(RQ)$ denotes the full embedding of the stalk subcategory at $\p$. Since the support of $\kp$ is equal to $\{\p\}$, the essential image of $\varphi_\p$ is contained in $\Gamma_\p\cD(RQ)$, here $\overline{\varphi}_\p$ denotes the corestriction. To justify the adjunction, first notice that $\varphi_\p=\iota_\p\circ\overline{\varphi}_\p$, then for any $X,Z\in\Gamma_\p\cD(RQ)$ and $Y\in\cD(\kp Q)$, we have isomorphisms
\[\Hom_{\Gamma_\p\cD(RQ)}(X,\overline{\varphi}_\p Y)\cong\Hom_{\cD(RQ)}(\iota_\p X,\iota_\p\overline{\varphi}_\p Y)\cong\Hom_{\cD(\kp Q)}(\varphi_\p^\ast\iota_\p X,Y)\]
and
\[\Hom_{\Gamma_\p\cD(RQ)}(\overline{\varphi}_\p Y,Z)\cong\Hom_{\cD(RQ)}(\iota_\p\overline{\varphi}_\p Y,\iota_\p Z)\cong\Hom_{\cD(\kp Q)}(Y,\varphi_\p^!\iota_\p Z).\]
By abuse of notation, we will denote this adjoint triple using the same notation as the one in the previous section.

Note that to study the closure property of t-structures in $\Gamma_\p\cD(RQ)$ under the composite $\varphi_\p\varphi_\p^\ast$ and
$\varphi_\p\varphi_\p^!$, we can not directly use the same arguments of \cite[Proposition 1.11]{Sab}. However we still get an analogous result.

\begin{lem}\label{Closure}
    For any t-structure $(\cU,\cW)$ in $\Gamma_\p\cD(RQ)$, it holds that:
    \begin{enumerate}
        \item $\varphi_\p\varphi_\p^\ast\cU=\kp\otimesL\cU\subseteq\cU$;
        \item $\varphi_\p\varphi_\p^!\cW=\RHom_R(\kp,\cW)\subseteq\cW$.
    \end{enumerate}
\begin{proof}
    (1) Recall that $\cU$ is a suspended subcategory of $\Gamma_\p\cD(RQ)$ closed under coproducts. Since $\Gamma_\p\cD(RQ)$ is a localizing subcategory of $\cD(RQ)$, $\cU$ is still suspended and closed under coproducts in $\cD(RQ)$. Thus, by \cite[Lemma 1.10 (1)]{Sab}, $\kp\otimesL\cU\subseteq\cU$.

    \noindent (2) We have $\varphi_\p\varphi_\p^!\cW\subseteq\cW$ if and only if $\Hom_{\cD(RQ)}({}^\perp\cW,\varphi_\p\varphi_\p^!\cW)=0$. By adjuction, this holds if and only if $\Hom_{\cD(RQ)}(\varphi_\p\varphi_\p^\ast{}^\perp\cW,\cW)=0$, which holds by point (1).
\end{proof}
\end{lem}

Now we want to prove that, for any prime $\p\in\Spec(R)$, the assignment of \cref{Injection} restricts to a bijection between definable coaisles of the stalk subcategory $\Gamma_\p(\cD(RQ))$ and the ones of $\cD(\kp Q)$. Recall form \cref{Gammap}, that stalk subcategories are compactly generated algebraic triangulated categories, thus the equivalence between homotopically smashing coaisles and cosuspended definable subcategories, mentioned in \cref{EqProp}, also holds in $\Gamma_\p(\cD(RQ))$.

\begin{thm}\label{Bij}
    Let $R$ be a commutative noetherian ring and $Q$ a finite quiver. Then, for any prime ideal $\p\in\Spec(R)$, there is an order preserving bijection
    \begin{equation*}\begin{gathered}
        F_\p:\Coaisle_{\Def}(\Gamma_\p\cD(RQ))\longrightarrow\Coaisle_{\Def}(\cD(\kp Q)) \\
        \cW=\cJ^\perp\mapsto\add_{\kp Q}\leftangle\varphi_\p^!\cW\rightangle=\left(\varphi_\p^\ast\cJ\right)^\perp
    \end{gathered}\end{equation*}
    with inverse $G_\p:\mathsf{w}\mapsto\cosusp_{RQ}^\hs\leftangle\varphi_\p \mathsf{w}\rightangle$. Moreover, any definable coaisle in $\Gamma_\p\cD(RQ)$ is compactly generated, in particular, the stalk subcategories $\Gamma_\p\cD(RQ)$ satisfy the telescope conjecture.
\begin{proof}
    The assignment $F_\p$ is well defined and injective. Indeed, for a definable coaisle $\cW=\cJ^\perp$ of $\Gamma_\p\cD(RQ)$, we can prove that $\add_{\kp Q}\leftangle\varphi_\p^!\cW\rightangle=\left(\varphi_\p^\ast\cJ\right)^\perp$ and $\cW=G_\p\circ F_\p(\cW)$ following the proof of \cref{Injection}, which still holds in this context by \cref{Closure}.

    \noindent Let us prove that $F_\p$ is surjective. Let $\mathsf{w}$ be a definable coaisle of $\cD(\kp Q)$, by the telescope conjecture result of \cite[Theorem 3.11]{AH} and its characterization in \cref{TCideal}, there is a set $\mathsf{s}$ of identity maps between compact objects of $\cD(\kp Q)$ such that $\mathsf{w}=\mathsf{s}^\perp$. Denote by $\cS=\leftcurly\id_{\widetilde{S}}\mid\id_S\in\mathsf{s}\rightcurly$, the set of identity maps between the lifts $\widetilde{S}\in(\Gamma_\p\cD(RQ))^c$ (see \cref{Lift}) for every $S$ such that $\id_S\in\mathsf{s}$. Since $\mathsf{s}$ is closed under positive shifts, $\cS^\perp$ is a (definable) coaisle of $\Gamma_\p\cD(RQ)$ and, by \cref{Lift} (2), $\mathsf{s}^\perp=(\varphi_\p^*\cS)^\perp$, i.e.~$\mathsf{w}=F_\p(\cS^\perp)$.

    \noindent As for the last part, it follows from the bijection that any definable coaisle $\cW$ of $\Gamma_\p\cD(RQ)$ is of the form $\cW=G_\p(\mathsf{s}^\perp)$, for some set $\mathsf{s}$ of identity maps between compact objects of $\cD(\kp Q)$. Moreover, $G_\p(\mathsf{s}^\perp)=G_\p\circ F_\p(\cS^\perp)=\cS^\perp$, for a set $\cS$ of identity maps between compact objects of $\Gamma_\p\cD(RQ)$, and thus it is compactly generated.
\end{proof}
\end{thm}

Let us now restrict the discussion to commutative artinian rings. In particular, for such a ring $R$, the prime spectrum is finite and discrete, consisting of finitely many maximal ideals $\Spec(R)=\{\m_1,\ldots,\m_n\}$ and the ring decomposes into the direct product of its localizations $R\cong\bigoplus_{i=1}^n R_{\m_i}$ (see \cite[Theorems 8.3, 8.5, 8.7]{AM}).

\begin{cor}\label{Main}
    For any commutative artinian ring $R$ and finite quiver $Q$, the derived category $\cD(RQ)$ satisfies the telescope conjecture.
\begin{proof}
    Let $R\cong\bigoplus_{i=1}^n R_{\m_i}$ be a commutative artinian ring. Note that, for each $\m\in\Spec(R)$, $R_\m$ is a 0-dimensional local ring and thus every $R_\m Q$-module is supported on $\{\m\}$, i.e.~$\cD(R_\m Q)\cong\Gamma_\m\cD(RQ)$ (see \cref{SuppProp}). The decomposition $R\cong\bigoplus_{i=1}^n R_{\m_i}$ of the ring induces a decomposition $\Mods_R(Q) \cong \Mods_{R_{\m_1}}(Q) \times\ldots\times\Mods_{R_{\m_n}}(Q)$ of the module category. It follows that the derived category $\cD(RQ)$ is equivalent to the product category $\Gamma_{\m_1}\cD(RQ)\times\ldots\times\Gamma_{\m_n}\cD(RQ)$ and so, by \cref{Bij}, it satisfies the telescope conjecture.
\end{proof}
\end{cor}

\begin{exmp}[Representations over truncated polynomial algebras]\label{Example}\hfill

\noindent An application of \cref{Main} to the representation theory of finite-dimensional algebras is provided by representations of finite acyclic quivers $Q$ over truncated polynomial algebras
\[\Lambda:=\bK[x_1,\ldots,x_m]/\left(x_i^{n_i}\mid 1\leq i\leq m\right)\]
where $\bK$ is a field and $n_i\geq 2$ for all $1\leq i\leq m$.

In this case $\Lambda Q$ is again a finite-dimensional algebra, and hence isomorphic to a path algebra $\bK Q'/I$ for some quiver $Q'$ and admissible ideal $I$, where
\[Q'_0=Q_0,\ Q'_1=Q_1\cup\leftcurly l_{1,k},\ldots,l_{m,k}\mid k\in Q'_0\rightcurly\text{ and}\]
\[I=\left(\alpha l_{i,k} - l_{i,h}\alpha,\,l_{i,k}^{n_i}\mid 1\leq i\leq m,\,k\in Q'_0,\,\alpha:k\to h\in Q'_1\right)\]
with $l_{i,k}$ denoting a loop at the vertex $k$ for any $1\leq i\leq m$. Therefore, these algebras belong among the few known examples of non-hereditary finite dimensional algebras satisfying the telescope conjecture.

For the so-called algebra of dual numbers $\Lambda=\bK[x]/(x^2)$, as discussed in \cite{RZ}, representations of a finite quiver $Q$ over $\Lambda$ coincides with the category of differential $\bK Q$-modules $\mathrm{Diff}(\bK Q)$. Thus, in this framework, \cref{Main} proves the telescope conjecture for the derived category $\cD(\mathrm{Diff}(\bK Q))$, or equivalently, for the $J$-shaped derived category $\cD_J(\bK Q)$ (see for example \cite[Section 4]{HJ}), where $J$ denotes the Jordan quiver \begin{tikzcd} \overset{1}{\bullet} \arrow[in=330, out=30, loop, "\delta", near start] \end{tikzcd} with the relation $\delta^2=0$.
\end{exmp}

\section{Towards non-noetherianity}\label{Sec4}
We want to dedicate this section to extending the previous results out of the noetherian world. To do this, we focus our attention on commutative rings which decompose into a product of $0$-dimensional local rings which satisfy a minimality condition as \cref{LTGM} (2), the \emph{perfect rings}. Recall that a commutative ring $R$ is \emph{semi-local} if it has only finitely many maximal ideals, and \emph{semi-artinian} if it admits a (possibly infinite) composition series -- that is, if $R$ as a regular $R$-module can be obtained as a transfinite extension of simple $R$-modules. Note that a semi-artinian ring is artinian if and only if it is noetherian. Moreover, an ideal $I$ of $R$ is said to be \emph{$T$-nilpotent} if, for every sequence $r_1,r_2,\ldots$ in $I$, there exists an $n>0$ such that $r_1r_2\ldots r_n=0$.

\begin{defn}[{\cite[Proposition 2.6]{BG}}]\label{PerfDef}
    A commutative ring $R$, is called \emph{perfect} if one of the following equivalent properties is satisfied:
    \begin{itemize}
        \item[(i)\,] $R$ is semi-local and semi-artinian;
        \item[(ii)] $R$ is a finite direct product of local rings with $T$-nilpotent maximal ideals;
        \item[(iii)] Every $R$-module has a projective cover.
    \end{itemize}
\end{defn}

In particular, if $R$ is a commutative perfect ring then, due to its semi-artinianity, it is $0$-dimensional (see \cite[Example 2.7]{KT}) and, due to its semi-local condition, one has $\Spec(R)=\{\m_1,\ldots,\m_n\}$ where the $\m_i$ are the maximal ideals of $R$. Moreover, by condition (ii), it follows that a commutative perfect ring is isomorphic to a finite direct product of local perfect rings each of them being the localization at a maximal ideal, i.e.~$R\cong\bigoplus_{i=1}^n R_{\m_i}$. Thus, for now, let us focus on local perfect rings.

Let us prove that for a local perfect ring $(R,\m,\kappa)$ the algebra $RQ$ satisfies a minimality condition for homotopically smashing complete cosuspended subcategories as the one stated in \cref{LTGM} (2).

\begin{prop}\label{PerfMin}
    Let $(R,\m,\kappa)$ be a commutative local perfect ring and $Q$ be a finite quiver. Then, for any complex $Y\in\cD(RC)$, it holds that
    \[\cosusp_{RQ}^{\Pi,\hs}\leftangle Y\rightangle=\cosusp_{RQ}^{\Pi,\hs}\leftangle\RHom_R(\kappa,Y)\rightangle.\]
\begin{proof}
    The proof follows as the one of \cite[Theorem 2.3]{Sab}. Indeed, by semi-artinianity, the ring $R$ itself is a transfinite extension of copies of the residue field $\kappa$. So, since $\kappa$ is contained in the cocomplete suspended subcategory
    \[\cM=\leftcurly M\in\cD(R)\mid\RHom_R(M, Y)\in\cosusp_{RQ}^{\Pi,\hs}\leftangle\RHom_R(\kappa, Y)\rightangle\rightcurly\]
    thus contains also $R$, and we can conclude that $Y \cong \RHom_R(R,Y)\in\cosusp_{RQ}^{\Pi,\hs}\leftangle\RHom_R(\kappa, Y)\rightangle$.
\end{proof}
\end{prop}

We will see that the bijection of \cref{Bij}, holds also for a local perfect ring $R$ but before doing that let us briefly recall some results of the previous sections in the present context.

\begin{prop}\label{Adapt}
    Let $(R,\m,\kappa)$ be a commutative local perfect ring, then:
    \begin{enumerate}
        \item The adjoint triple $\varphi_\m^\ast\dashv\varphi_\m\dashv\varphi_\m^!$ satisfies:
            \begin{enumerate}
                \item $\varphi_\m\circ\varphi_\m^!\cong\RHom_R(\kappa,\_)$;
                \item $\varphi_\m^!\circ\varphi_\m(X)\cong\prod_{i\geq0}\kappa^{(\alpha_i)}[i]$ for some cardinals $\alpha_i$.
            \end{enumerate}
        \item For any compact object $S\in\cD(\kp Q)$ there exists $\widehat{S}\in\cD^c(RQ)$ such that $\kappa\otimes_R\widehat{S}\cong S$.
    \end{enumerate}
\begin{proof}
    (1) Point (a) follows as in \cref{AdjTrip} (1). Since $R$ is local and perfect, by condition (iii) of \cref{PerfDef}, $\kappa$ admits a minimal free resolution $\mathbf{f}\kappa$ and $\kappa\otimes_R\mathbf{f}\kappa\cong\bigoplus_{i\geq0}\kappa^{(\alpha_i)}[i]$ for some cardinals $\alpha_i$. Thus, (b) follows as \cref{AdjTrip} (2) with the obvious modifications.

    \noindent(2) The result follows as \cref{Lift} (1).
\end{proof}
\end{prop}

\begin{thm}\label{LocPerf}
    Let $(R,\m,\kappa)$ be a commutative local perfect ring and $Q$ a finite quiver. Then, there is an order preserving bijection
    \begin{equation*}\begin{gathered}
        F:\Coaisle_{\Def}(\cD(RQ))\longrightarrow\Coaisle_{\Def}(\cD(\kappa Q))\\
        \cW=\cJ^\perp\mapsto\add_{\kappa Q}\leftangle\varphi_\m^!\cW\rightangle=\left(\varphi_\m^\ast\cJ\right)^\perp
    \end{gathered}\end{equation*}
    with inverse $G:\mathsf{w}\mapsto\cosusp_{RQ}^{\Pi,\hs}\leftangle\varphi_\m \mathsf{w}\rightangle$. Moreover, any definable coaisle in $\cD(RQ)$ is compactly generated.
\begin{proof}
    That the assignment $F$ is well-defined follows by the same proof of \cref{Injection} (1), using \cref{Adapt} (1). As for the injectivity, given a definable coaisle $\cW\subseteq\cD(RQ)$, the subcategory $G\circ F(\cW)=\cosusp_{RQ}^{\Pi,\hs} \leftangle\varphi_\m\varphi_\m^!\cW\rightangle$ is contained in $\cW$ and, by \cref{PerfMin}, also the reverse containment holds. 

    \noindent The surjectivity of $F$ follows similarly to \cref{Bij}. Given a definable coaisle $\mathsf{w}$ of $\cD(\kappa Q)$, by the telescope conjecture result of \cite[Theorem 3.11]{AH} and \cref{TCideal}, there is a set $\mathsf{s}\subseteq\cD^c(\kappa Q)$ of identity maps such that $\mathsf{w}=\mathsf{s}^\perp$. Denote by $\cS=\leftcurly\id_{\widehat{S}}\mid\id_S\in\mathsf{s}\rightcurly$ the set of identity maps between the lifts $\widehat{S}\in\cD^c(RQ))$ for every $S$ such that $\id_S\in\mathsf{s}$. Since $\mathsf{s}$ is closed under positive shifts, $\cS^\perp$ is a (definable) coaisle of $\cD(RQ)$ and, by \cref{Adapt} (2), $\mathsf{s}^\perp=(\varphi_\p^*\cS)^\perp$, i.e.~$\mathsf{w}=F_\p(\cS^\perp)$.

    \noindent As for the last part, it follows from the bijection that any definable coaisle $\cW$ of $\cD(RQ)$ is of the form $\cW=G(\mathsf{s}^\perp)$, for some set of identity maps $\mathsf{s}\subseteq\cD^c(\kappa Q)$. Moreover, $G(\mathsf{s}^\perp)=G\circ F(\cS^\perp)=\cS^\perp$, for a set of identity maps $\cS\subseteq\cD^c(RQ)$, and thus it is compactly generated.
\end{proof}
\end{thm}

So, analogously to \cref{Main}, we can conclude the following.

\begin{cor}
    For any commutative perfect ring $R$ and finite quiver $Q$, the derived category $\cD(RQ)$ satisfies the telescope conjecture.
\begin{proof}
    Let $R$ be a commutative perfect ring, then by definition it is equivalent to a direct product of finitely many local perfect rings $R\cong\bigoplus_{i=1}^nR_i$. It follows again that the derived category $\cD(RQ)$ is equivalent to the product category $\cD(R_1 Q)\times\ldots\times\cD(R_n Q)$ and so, by \cref{LocPerf}, it satisfies the telescope conjecture.
\end{proof}
\end{cor}

\raggedbottom
\bibliographystyle{amsalpha}
\bibliography{references}

@book{Pop,
    author={Popescu, N.},
    title={Abelian categories with applications to rings and modules},
    publisher={Academic Press},
    year={1973},
}

@book{CFH,
    author={Christensen, L.~W. and Foxby, H.~B. and Holm, H.},
    title={Derived category methods in commutative algebra},
    publisher={Springer Monographs in Mathematics}, 
    year={2024},
}

@book{Hov,
    author={Hovey, M.},
    title={Model Categories},
    series={Mathematical Surveys and Monographs},
    volume={63},
    publisher={American Mathematical Society},
    year={1999},
}

@book{Ben,
    author={Benson, D.~J.},
    title={Representations and cohomology, vol.~{I}},
    publisher={Cambridge University Press}, 
    year={1995},
}

@book{AM,
    author={Atiyah, M.~F. and MacDonald, I.~G.},
    title={Introduction to commutative algebra},
    publisher={CRC Press},
    year={1969},
}

@misc{Sab,
    author={Sabatini, E.},
    title={Telescope conjecture for t-structures over noetherian path algebras},
    howpublished={arXiv preprint \href{https://arxiv.org/abs/2505.20803}{arXiv:2505.20803}},
    year={2025},
}

@article{Bou,
    author={Bousfield, A.~K.},
    title={The localization of spectra with respect to homology},
    journal={Topology},
    volume={18},
    number={4},
    pages={257--281},
    year={1979},
}

@article{Rav,
    author={Ravenel, D.~C.},
    title={Localization with respect to certain periodic homology theories},
    journal={American Journal of Mathematics},
    volume={106},
    number={2},
    pages={351--414},
    year={1984},
}

@article{NeeTC,
    author={Neeman, A.},
    title={The chromatic tower for {D(R)} \rm{(with an appendix by Marcel B{\"o}kstedt)}},
    journal={Topology},
    volume={31},
    number={3},
    pages={519--532},
    year={1992},
}

@article{BHTC,
    author={Bazzoni, S. and Hrbek, M.},
    title={Definable coaisles over rings of weak global dimension at most one},
    journal={Publicacions Matem{\`a}tiques},
    volume={65},
    number={1},
    pages={165--241},
    year={2021},
}

@article{KS,
    author={Krause, H. and {\v S}{\v t}ov{\'i}{\v c}ek, J.},
    title={The telescope conjecture for hereditary rings via {E}xt-orthogonal pairs},
    journal={Advances in Mathematics},
    volume={225},
    number={5},
    pages={2341--2364},
    year={2010},
}

@article{Spa,
    title={Resolutions of unbounded complexes},
    author={Spaltenstein, N.},
    journal={Compositio Mathematica},
    volume={65},
    number={2},
    pages={121--154},
    year={1988},
}

@article{AS,
    author={Antieau, B. and Stevenson, G.},
    title={Derived categories of representations of small categories over commutative noetherian rings},
    journal={Pacific Journal of Mathematics},
    volume={283},
    number={1},
    pages={21--42},
    year={2016},
}

@article{Nee,
    author={Neeman, A.},
    title={The {G}rothendieck duality theorem via {B}ousfield’s techniques and {B}rown representability},
    journal={Journal of the American Mathematical Society},
    volume={9},
    number={1},
    pages={205--236},
    year={1996},
}

@article{HN,
    author={Hrbek, M. and Nakamura, T.},
    title={Telescope conjecture for homotopically smashing t-structures over commutative noetherian rings},
    journal={Journal of Pure and Applied Algebra},
    volume={225},
    number={4},
    pages={106571},
    year={2021},
}

@incollection{KraAlg,
    author={Krause, H.},
    title={Derived categories, resolutions, and {B}rown representability},
    booktitle={Interactions between homotopy theory and algebra},
    series={Contemporary Mathematics},
    volume ={436},
    publisher={American Mathematical Society},
    pages={101--139},
    year={2007},
}

@incollection{Kel,
    author={Keller, B.},
    title={On differential graded categories},
    booktitle={International Congress of Mathematicians (Madrid 2006), vol.~{II}},
    pages={151--190},
    publisher={European Mathematical Society},
    year={2006},
}

@article{Rin,
    author={Ringel, C.~M.},
    title={Exceptional modules are tree modules},
    journal={Linear Algebra and its Applications},
    volume={275-276},
    pages={471--493},
    year={1998},
}

@article{CB,
    author={Crawley-Boevey, W.},
    title={Rigid integral representations of quivers over arbitrary commutative rings}, 
    journal={Annals of Representation Theory},
    volume={1},
    number={3},
    pages={375--384},
    year={2024},
}

@article{Bec,
    author={Becker, H.},
    title={Models for singularity categories},
    journal={Advances in Mathematics},
    volume={254},
    pages={187--232},
    year={2014},
}

@article{SSV,
    author={Saor{\'i}n, M. and {\v S}{\v t}ov{\'i}{\v c}ek, J. and Virili, S.},
    title={t-{S}tructures on stable derivators and {G}rothendieck hearts},
    journal={Advances in Mathematics},
    volume={429},
    pages={109139},
    year={2023},
}

@article{KV,
    author={Keller, B. and Vossieck, D.},
    title={Aisles in derived categories},
    journal={Bulletin of the Belgian Mathematical Society},
    volume={40},
    number={2},
    pages={239--253},
    year={1988},
}

@article{NeeSet,
    author={Neeman, A.},
    title={The t-structures generated by objects},
    journal={Transactions of the American Mathematical Society},
    volume={374},
    number={11},
    pages={8161--8175},
    year={2021},
}

@article{KraSma,
    author={Krause, H.},
    title={Smashing subcategories and the telescope conjecture -- an algebraic approach},
    journal={Inventiones Mathematicae},
    volume={139},
    pages={99--133},
    year={2000},
}

@article{SS,
    author={Saor{\'i}n, M. and {\v S}{\v t}ov{\'i}{\v c}ek, J.},
    title={t-{S}tructures with {G}rothendieck hearts via functor categories},
    journal={Selecta Mathematica},
    volume={29},
    number={5},
    pages={77},
    year={2023},
}

@article{Kra,
    author={Krause, H.},
    title={Cohomological quotients and smashing localizations},
    journal={American Journal of Mathematics},
    volume={127},
    number={6},
    pages={1191--1246},
    year={2005},
}

@article{LV,
    author={Laking, R. and Vit{\'o}ria, J.},
    title={Definability and approximations in triangulated categories},
    journal={Pacific Journal of Mathematics},
    volume={306},
    number={2},
    pages={557--586},
    year={2020},
}

@article{CS,
    title={Uniqueness of dg enhancements for the derived category of a {G}rothendieck category},
    author={Canonaco, A. and Stellari, P.},
    journal={Journal of the European Mathematical Society},
    volume={20},
    number={11},
    pages={2607--2641},
    year={2018},
}

@incollection{Stour,
    author={Stevenson, G.},
    title={A Tour of support theory for triangulated categories through tensor triangular geometry},
    booktitle={Building Bridges Between Algebra and Topology},
    series={Advanced Courses in Mathematics - CRM Barcelona},
    publisher={Birkh{\"a}user},
    pages={63--101},
    year={2018},
}

@article{AH,
    author={Angeleri{ }H{\"u}gel, L. and Hrbek, M.},
    title={Parametrizing torsion pairs in derived categories},
    journal={Representation Theory of the American Mathematical Society},
    volume={25},
    number={23},
    pages={679--731},
    year={2021},
}

@article{RZ,
    title={Representations of quivers over the algebra of dual numbers},
    author={Ringel, C.~M. and Zhang, P.},
    journal={Journal of Algebra},
    volume={475},
    pages={327--360},
    year={2017},
}

@article{HJ,
    title={Minimal semiinjective resolutions in the {Q}-shaped derived category},
    author={Holm, H. and J{\o}rgensen, P.},
    journal={Revista Matem{\'a}tica Iberoamericana},
    pages={1--26},
    year={2025},
}

@article{BG,
    title={Covering classes and 1-tilting cotorsion pairs over commutative rings},
    author={Bazzoni, S. and Gros, G.~Le},
    journal={Forum Mathematicum},
    volume={33},
    number={3},
    pages={601--629},
    year={2021},
}

@article{KT,
    author={Kourki, F. and Tribak, R.},
    title={On semiartinian and {$\Pi$}-semiartinian modules},
    journal={Palestine Journal of Mathematics},
    volume={7},
    pages={99--107},
    year={2018},
}
\end{document}